\documentclass[11pt]{article}

\usepackage{amsmath, amssymb, amscd}

\def\eqref#1{(\ref{#1})}
\newcommand{\goth}{\frak}

\newcommand{\arrow}{{\:\longrightarrow\:}}
\newcommand{\Z}{{\Bbb Z}}
\newcommand{\C}{{\Bbb C}}
\newcommand{\R}{{\Bbb R}}
\newcommand{\Q}{{\Bbb Q}}

\newcommand{\6}{\partial}
\def\1{\sqrt{-1}\:}
\newcommand{\restrict}[1]{{\left|_{{\phantom{|}\!\!}_{#1}}\right.}}

\newcommand{\calo}{{\cal O}}

\renewcommand{\bar}{\overline}
\renewcommand{\phi}{\varphi}
\renewcommand{\epsilon}{\varepsilon}
\renewcommand{\geq}{\geqslant}
\renewcommand{\leq}{\leqslant}

\newcommand{\End}{\operatorname{End}}
\newcommand{\Tot}{\operatorname{Tot}}

\newcommand{\Aut}{\operatorname{Aut}}
\newcommand{\Iso}{\operatorname{Iso}}

\newcommand{\Diff}{\operatorname{Diff}}

\newcommand{\Lie}{\operatorname{Lie}}

\newcommand{\comment}[1]{{}}

\def\blacksquare{\hbox{\vrule width 4pt height 4pt depth 0pt}}
\def\endproof{\blacksquare}

\makeatletter

\@ifundefined{Bbb}
     {\newcommand{\Bbb}[1]{{\mathbb #1}}}%
{}%     {\edef\Bbb#1{{\Bbb #1}}}

\newcommand{\ps@verbit}{%
  \renewcommand{\@oddhead}{%
          \scriptsize
          {Immersion theorem for Vaisman manifolds}
          \hfil\tiny {L. Ornea and M. Verbitsky, May 30, 2003 }}
  \renewcommand{\@evenhead}{\@oddhead}
  \renewcommand{\@oddfoot}{\hfil\thepage\hfil}
  \renewcommand{\@evenfoot}{\@oddfoot}}
 
\newcounter{Mycounter}[section]
\newcounter{lemma}[section]
\renewcommand{\thelemma}{\noindent{Lemma \thesection.\arabic{lemma}}}
\newcommand{\lemma}{%
     \setcounter{lemma}{\value{Mycounter}}
     \refstepcounter{lemma}
     \stepcounter{Mycounter}
     {\bf \thelemma:\ }}

\newcounter{claim}[section]
\renewcommand{\theclaim}{\noindent{Claim \thesection.\arabic{claim}}}
\newcommand{\claim}{%
     \setcounter{claim}{\value{Mycounter}}
     \refstepcounter{claim}
     \stepcounter{Mycounter}
     {\bf \theclaim:\ }}

\newcounter{sublemma}[section]
\newcounter{corollary}[section]
\renewcommand{\thecorollary}{\noindent{Corollary
 \thesection.\arabic{corollary}}}
\newcommand{\corollary}{%
     \setcounter{corollary}{\value{Mycounter}}
     \refstepcounter{corollary}
     \stepcounter{Mycounter}
     {\bf \thecorollary:\ }}

\newcounter{theorem}[section]
\renewcommand{\thetheorem}{\noindent{Theorem
 \thesection.\arabic{theorem}}}
\newcommand{\theorem}{%
     \setcounter{theorem}{\value{Mycounter}}
     \refstepcounter{theorem}
     \stepcounter{Mycounter}
     {\bf \thetheorem:\ }}

\newcounter{conjecture}[section]
\newcounter{proposition}[section]
\renewcommand{\theproposition}
       {\noindent{Proposition \thesection.\arabic{proposition}}}
\newcommand{\proposition}{%
     \setcounter{proposition}{\value{Mycounter}}
     \refstepcounter{proposition}
     \stepcounter{Mycounter}
     {\bf \theproposition:\ }}

\newcounter{definition}[section]
\renewcommand{\thedefinition}
       {\noindent{Definition~\thesection.\arabic{definition}}}
\newcommand{\definition}{%
     \setcounter{definition}{\value{Mycounter}}
     \refstepcounter{definition}
     \stepcounter{Mycounter}
     {\bf \thedefinition:\ }}

\newcounter{example}[section]
\newcounter{remark}[section]
\renewcommand{\theremark}{\noindent{Remark \thesection.\arabic{remark}}}
\newcommand{\remark}{%
     \setcounter{remark}{\value{Mycounter}}
     \refstepcounter{remark}
     \stepcounter{Mycounter}
     {\bf \theremark:\ }}

\newcounter{problem}[section]
\newcounter{question}[section]
\renewcommand{\thequestion}{\noindent{Question
 \thesection.\arabic{question}}}
\newcommand{\question}{%
     \setcounter{question}{\value{Mycounter}}
     \refstepcounter{question}
     \stepcounter{Mycounter}
     {\bf \thequestion:\ }}

\@addtoreset{equation}{section}
\@addtoreset{footnote}{section}
\makeatother

\begin{document}

%%%%%%%%%%%%%%%%%%%%%%%%%%%%%%%%%%%%%%%%%%%%%%%%%%%%%%%%%%%%
\begin{center}
{\LARGE\bf
Immersion theorem for Vaisman manifolds 
}
%%%%%%%%%%%%%%%%%%%%%%%%%%%%%%%%%%%%%%%%%%%%%%%%%%%%%%%%%%%%
\\[4mm]
Liviu Ornea\footnote{Liviu Ornea is member of EDGE, Research
Training Network HRPN-CT-2000-00101, supported by the European Human
Potential Programme.}, Misha Verbitsky,\footnote{Misha Verbitsky is 
an EPSRC advanced fellow 
supported by CRDF grant RM1-2354-MO02 and EPSRC grant  GR/R77773/01

Both authors acknowledge financial support from Ecole Polytechnique
 (Palaiseau). \\[1mm]

2000 {\it MSC}: 53C55, 14E25, 53C25.

}
\\[4mm]

%{\tt lornea@imar.ro; \ \  verbit@maths.gla.ac.uk, \ \  verbit@mccme.ru}
\end{center}

%%%%%%%%%%%%%%%%%%%%%%%%%%%%%%%%%%%%%%%%%%%%%%%%
{\small 
\hspace{0.15\linewidth}
\begin{minipage}[t]{0.7\linewidth}
{\bf Abstract} \\
A locally conformally K{\"a}hler 
(LCK) manifold is a complex manifold
admitting a K{\"a}hler covering $\tilde M$, 
with monodromy acting on $\tilde M$ by K{\"a}hler homotheties.
A compact LCK manifold is Vaisman if it admits a holomorphic
flow acting by non-trivial homotheties on $\tilde M$.
We prove that any compact Vaisman manifold  
admits a natural holomorphic immersion to a 
Hopf manifold $(\C^n\setminus 0)/\Z$.
As an application, we obtain that any Sasakian
manifold has a contact immersion to an 
odd-dimensional sphere.
\end{minipage}
}
%%%%%%%%%%%%%%%%%%%%%%%%%%%%%%%%%%%%%%%%%%%%%%%%

{
\small
\tableofcontents
}
%%%%%%%%%%%%%%%%%%%%%%%%%%%%%%%%%%%%%%%%%%%%%%%%

\section{Introduction}
\label{_Intro_Section_}

\subsection{LCK manifolds: a historical overview}

Locally conformally K{\"a}hler (LCK)
manifolds are, by definition,
complex manifolds admitting a K{\"a}hler covering  
with deck transformations acting by K{\"a}hler homotheties. 
These manifolds appear naturally in complex geometry.
Most examples of compact non-K{\"a}hler manifolds studied
in complex geometry admit an LCK structure. 

LCK structures appeared since 1954, studied by P.
Libermann. They came again into attention since 1976, with the work of
I. Vaisman who also produced the first compact examples: the
(diagonal) Hopf manifolds (see \cite{_Dragomir_Ornea_} and the
references therein). 

As any Hermitian geometry, LCK geometry is encoded in the properties
of the Lee form (\ref{_LCK_Definition_})
which, in this case, is closed. From the conformal
viewpoint, the Lee form is the connection one-form in a real line
 bundle,
called the weight bundle. 

It became soon clear that the most tractable LCK manifolds are the
ones with parallel Lee form. We call them Vaisman manifolds because it
was I. Vaisman to introduce and study them extensively,
under the name of generalized Hopf manifolds. With the
exception of some  compact surfaces (Inoue surfaces of the first kind,
Hopf surfaces of K{\"a}hler rank 0) which are known to not admit Vaisman
metrics (cf. \cite{_Belgun_}), all other known examples of 
compact LCK manifolds
have parallel Lee form. 

The universal covering space of a Vaisman manifold bears a K{\"a}hler
metric of a very special kind: it is a conic metric over a Sasakian
manifold. Sasakian spaces proved recently to be very important
in physical theories (see the series of papers by C.P. Boyer,
K. Galicki \emph{et.al.}); this further motivated the interest for
Vaisman manifolds. 

The topology of (compact) Vaisman manifolds is very different from that
 of
K{\"a}hler manifolds. While it was long ago known that the first Betti
number must be odd, \cite{_Vaisman:Dedicata_}, it was proved only
recently that Einstein-Weyl LCK manifolds
have no non-trivial holomorphic forms, 
\cite{_Aleksandrov_Ivanov_}. 
Still, nothing is known about the topology 
of general LCK manifolds.  

The parallelism of the Lee form is 
reflected in the properties of the
weight bundle: indeed, we recently proved in
\cite{_OV:Structure_} that the weight bundle has
monodromy $\Z$. This allowed us to obtain the structure theorem of
compact Vaisman manifolds: they are Riemannian suspensions over the
circle with fibre a Sasakian manifold. This result can be regarded as
a kind of equivalence between Sasakian and Vaisman geometries.

Vaisman manifolds can be also viewed \emph{via} the properties of
their automorphism group. This was done separately by F.A. Belgun and
Y. Ka\-mi\-shi\-ma. Finally, the existence of a complex flow of
 conformal
transformations proved to be an equivalent definition of Vaisman
structures, cf. \cite{_Kamishima_Ornea_}. This was also the approach
in \cite{_Gini_Ornea_Parton_} where the Hamiltonian actions were
considered.

In the present paper we further exploit the techniques developed in
\cite{_Verbitsky:LCHK_}. We adapt algebraic geometrical methods to
prove an analogue of the Kodaira embedding theorem in Vaisman
geometry.

%%%%%%%%%%%%%%%%%%%%%%%%%%%%%%%%%%%%%%%%%%%%%%%%%%%%%%%%%%%%
\subsection{Algebraic geometry of LCK manifolds}
%%%%%%%%%%%%%%%%%%%%%%%%%%%%%%%%%%%%%%%%%%%%%%%%%%%%%%%%%%%%

The fundamental group $\pi_1(M)$ of an LCK manifold $M$
 acts on the K{\"a}hler covering $\tilde M$ by homotheties.
This gives a representation $\rho:\; \pi_1(M) \arrow \R^{>0}$.
{}From algebro-geometric point of view, the most interesting feature
of an LCK manifold is its weight bundle $L$, that is, the flat line
bundle with the monodromy defined by $\rho$.
This bundle is real, but its complexification $L_\C$
is flat, and therefore holomorphic. 

Since $L$ is a $\R^{>0}$-bundle, it is topologically trivial,
and admits a trivialization. The bundle $L$ can be interpreted
as a bundle of metrics in the conformal class of the LCK structure
on $M$. Each trivialization gives us a choice of a metric on $M$. 

Complex algebraic geometry deals mostly with the K{\"a}hler manifolds.
On a compact K{\"a}hler manifold $X$, flat line bundles are not
 particularly
interesting. They have no non-trivial holomorphic sections;
moreover, every flat bundle on $X$ is induced from the
Albanese map $X\arrow Alb(X)$ to the torus $Alb(X)$.

In non-K{\"a}hler complex geometry, the situation is 
completely different. The weight bundle $L_\C$
on a compact LCK manifold $M$ may admit quite a few
non-trivial sections. In some cases there are so many
sections that they provide an embedding from $M$
to a model LCK manifold, called a Hopf manifold
(\ref{_Hopf_mfld_Definition_}), in a manner of Kodaira
embedding theorem. 

To study the bundle $L_\C$, we fix a Hermitian metric
 (this is done using the Gauduchon's theorem,
see \ref{_Gaud_metric_Claim_}), and compute the
curvature $\Theta\in \Lambda^{1,1}(M)$
of the Chern connection on $L_\C$.
It turns out that all eigenvalues 
of $\Theta$ are positive except one.

J.-P. Demailly's
holomorphic Morse inequalities
(\cite{_Demailly:Morse_}) give
strong estimates of the asymptotic
behaviour of cohomology \[ H^i(L_\C^k), \ \  k \arrow \infty\]
for such bundles.

Even more interesting, the
weight covering $\tilde M$ of $M$ admits a
function $\psi:\; \tilde M \arrow \R^{>0}$
which is exhausting, and all  eigenvalues
of the form $-\1 \6\bar\6 \psi$ are positive,
except one. This function is obtained
by adding the K{\"a}hler potential $r= e^{-f}$ of $\tilde M$
(see \eqref{_kahler_pote_Equation_}) and $r^{-1}$.
Such a manifold is called 2-complete. There is a
great body of work dedicated to the study
of $q$-complete manifolds (see e.g.
the survey \cite{_Coltoiu_}).

The LCK manifolds are clearly an extremely interesting
object of al\-geb\-ro-\-geo\-met\-ric study.

%%%%%%%%%%%%%%%%%%%%%%%%%%%%%%%%%%%%%%%%%%%%%%%%
\subsection{Vaisman manifolds}
%%%%%%%%%%%%%%%%%%%%%%%%%%%%%%%%%%%%%%%%%%%%%%%%

In this paper, we deal with a special kind
of LCK manifolds, called Vaisman manifolds
(\ref{_Vaisman_Definition_}). 
Among other compact LCK manifolds, the Vaisman manifolds
are characterized as admitting a holomorphic flow
which acts by non-trivial homotheties on its
K{\"a}hler covering (see \cite{_Kamishima_Ornea_}). 

The weight monodromy of a compact Vaisman manifold
$M$ is isomorphic to $\Z$ (\cite{_OV:Structure_}).
Therefore, $M$ is a quotient $\tilde M/\Gamma$ of a 
K{\"a}hler manifold $\tilde M$ by $\Gamma \cong \Z$, where
$\Z$ acts on $\tilde M$ by holomorphic homotheties.

In \cite{_OV:Structure_}
it was proven that the K{\"a}hler covering $\tilde M$ of $M$ admits
an action of a commutative connected lie group
$\tilde G$, such that the monodromy $\Gamma$ lies
inside $\tilde G$.  We obtain that 
$\Gamma$ belongs to the connected component
of the group of holomorphic homotheties of $\tilde M$.
The converse is also true, by \cite{_Kamishima_Ornea_}:
if $\Gamma$ belongs to the connected component
of the group of holomorphic homotheties of $\tilde M$,
then $M$ is Vaisman. 

This gives the following characterization of
Vaisman manifolds. Given a K{\"a}hler manifold
$\tilde M$, denote by $H$ the connected
component of the Lie group of holomorphic
homotheties, and let $\Gamma\subset H$, $\Gamma\cong \Z$ 
be a group acting by non-trivial homotheties in such
a way that $M:= \tilde M/\Gamma$ is compact. Then
$M$ is Vaisman, and, moreover, all Vaisman
manifolds are obtained this way.

%%%%%%%%%%%%%%%%%%%%%%%%%%%%%%%%%%%%%%%%%%%%%%%%%%%%%%%%%%%%
\subsection{Quasiregular Vaisman manifolds}
\label{_qr_Vaisman_Intro_Subsection_}
%%%%%%%%%%%%%%%%%%%%%%%%%%%%%%%%%%%%%%%%%%%%%%%%%%%%%%%%%%%%

Let $M$ be a Vaisman manifold.
In \cite{_Verbitsky:LCHK_}, the curvature $\Theta$
of the Chern connection on the weight bundle was computed.
It was shown that $\Theta = -\1 \omega_0$, where 
$\omega_0$ is positive semidefinite form, with 
one zero eigenvalue. 

The zero eigenvalue of $\omega_0$ corresponds to a 
holomorphic foliation $\Xi$ on $M$, called {\bf the Lee foliation}. 
The leaves of this 
foliation are orbits of a holomorphic flow on $M$, called
{\bf the complex Lee flow} (see Subsection 
\ref{_qr_Vais_Subsection_}). 

When this foliation has a Hausdorff leaf space $Q$, one
would expect that the bundle $L_\C$  is a pullback
of a positive line bunlde $L_Q$ on $Q$. This is true 
(\ref{_L_trivi_on_fibe_Theorem_}), at least when the 
leaves of $\Xi$ are compact. 

Vaisman manifolds with compact leaves of the Lee foliation
are called {\bf quasiregular}. The manifold $Q$ is obtained
as a quotient of a manifold by an action of a compact group,
and therefore it is an orbifold. The form $\omega_0$ on $M$
defines a K{\"a}hler form on $Q= M/\Xi$. A compact orbifold
equipped with a positive line bundle is actually
projective, because Kodaira embedding theorem 
holds in the orbifold case (\cite{_Baily_}).

This approach allows for a completely algebraic construction
of quasiregular Vaisman manifolds. Consider a projective
orbifold $Q$ and a positive line bundle $L_Q$ on $Q$.
Let $\tilde M$ be the total space of all non-zero vectors
in $L_Q^*$, $q\in \R^{>0}$, and 
$\sigma_q:\; L_Q \arrow L_Q$ a map multiplying
$l \in L_Q\restrict{x\in M}$ by $q$.
Denote by $\Gamma \cong \Z$ the subgroup
of $\Aut(\tilde M)$ generated by $\sigma_q$.
When $\tilde M$ is smooth, it is K{\"a}hler,
because the function $\psi(b)= |b|^2$
gives a K{\"a}hler potential on 
$\tilde M$.\footnote{This is true 
(and easy to prove) for 
a space of non-zero vectors in
any negative bundle. For a proof, one may
check e.g. \cite{_Besse:Einst_Manifo_},
(15.19).} The quotient $M:= \tilde M/\sigma_q$
is clearly a Vaisman manifold; the Lee fibration
$M \arrow M/\Xi$ becomes the standard projection
to $Q$, and ths fibers of $L_Q$ are identified
with the leaves of Lee foliation. We obtain that
$M$ is quasiregular; moreover, all quasiregular
Vaisman manifolds are obtained this way.

In a sense, the quasiregular Vaisman geometry
is more algebraic than for instance K{\"a}hler
geometry. Indeed, any quasiregular Vaisman
manifold is constructed from a projective
manifold and an ample bundle - purely algebraic
set of data. 

Consider the lifting $\tilde L_\C^*$ of $L_\Q^*$ 
to \[ \tilde M \cong \Tot(L_Q^* \backslash 0).\]
The points of $\Tot(L_Q^* \backslash 0)$ correspond
to non-zero vectors in $L_\C^*$. Therefore,
the bundle  $\tilde L_\C^*$ is equipped
with a natural holomorphic non-degenerate 
section $\tau$.

Sections of $L_Q$ give, after a pullback, sections
of $\tilde L_\C$. Evaluating these sections
on $\tau$, we obtain holomorphic functions
on $\tilde M$. The map $\sigma_q$ multiplies
$\tau$ by $q$. Therefore, sections of $L_Q$
give functions $\mu$ on $\tilde M$ which satisfy
$\sigma_q(\mu) = q\mu$. Similarly, for sections
of $L_Q^k$, we have 
\begin{equation}\label{_sigma_q_on_L_Equation_}
\sigma_q(\mu) = q^k\mu.
\end{equation}
The bundle $L_Q$ is ample, and therefore,
for $k$ sufficiently big, Kodaira embedding theorem
gives an embedding $Q \hookrightarrow {\Bbb P}(H^0(L_Q^k))$.
Associating to sections of $L_Q^k$ functions on $\tilde M$
as above, we obtain a morphism
\begin{equation}\label{_tilde_M_immer_Equation_}
\tilde M \arrow H^0(L_Q^k)\backslash 0
\end{equation}
which maps the fibers of the Lee foliation to
complex lines of $H^0(L_Q^k)$ and induces 
an immersion. The automorphism
$\sigma_q$ acts on both sides of
\eqref{_tilde_M_immer_Equation_} as prescribed
by \eqref{_sigma_q_on_L_Equation_}, and this
gives an immersion
$M \arrow (H^0(L_Q^k)\backslash 0 )/\langle q^k \rangle$, 
where $\langle q^k \rangle\cong \Z$ denotes the abelian
subgroup of $\Aut (H^0(L_Q^k))$ generated by $q^k$.

\hfill

We outlined the proof of the following theorem.

\hfill

%%%%%%%%%%%%%%%%%%%%%%%%%%%%%%%%%%%%%%%%%%%%%%%%
\theorem \label{_qr_imme_Intro_Theorem_}
Let $M$ be a quasiregular Vaisman manifold. Then
$M$ admits a holomorphic immersion to $\C^n \backslash 0 / q$,
$q\in \R^{>1}$.

\hfill

For a detailed proof of \ref{_qr_imme_Intro_Theorem_},
see Subsection \ref{_qr_imme_Subsection_}.

\hfill

A similar theorem is true for general Vaisman manifolds
(\ref{_imme_the_gene_Theorem_}). 

%%%%%%%%%%%%%%%%%%%%%%%%%%%%%%%%%%%%%%%%%%%%%%%%%%%%%%%%%%%%
\subsection{Immersion theorem for general Vaisman manifolds}
%%%%%%%%%%%%%%%%%%%%%%%%%%%%%%%%%%%%%%%%%%%%%%%%%%%%%%%%%%%%

To deal with a Vaisman manifold $M$ which is not quasiregular,
we look at the Lie group $G$ within the group of isometries 
of $M$ generated by the Lee flow (see \ref{_Vaisman_Subsection_}).
This group was studied at great length in 
\cite{_Kamishima_} and
\cite{_Kamishima_Ornea_}.

It is not very difficult to show that
$G$ is a compact abelian group
(\cite{_Kamishima_Ornea_}). 
Denote by $\tilde G$ the group
of automorphisms of the pair $(\tilde M, M)$
and mapping to $G$ under the natural forgetful
map $\Aut(\tilde M, M) \arrow \Aut(M)$ (see 
Subsection \ref{_tilde_G_intro_Subsection_}). 
It was shown in \cite{_OV:Structure_}
that $\tilde G$ is connected and isomorphic
to $\R \times (S^1)^k$, and the natural forgetful map
$\tilde G \arrow G$ is a covering.

The deck transformation group  of $\tilde M$ (denoted by 
$\Gamma$) lies in $\tilde G$ (\cite{_OV:Structure_}; see also Subsection
\ref{_tilde_G_intro_Subsection_}).

For all sufficiently small
deformations $\Gamma'$ of $\Gamma$ within a 
complexification $\tilde G_\C \subset \Aut(\tilde M, M)$,
the quotient manifold $\tilde M / \Gamma'$ 
remains Vaisman (\ref{_defo_LCK_Theorem_}). 
Moreover, one can chose $\Gamma'$ in such a way
that $\tilde M/\Gamma'$ becomes quasiregular
(\ref{_qr_defo_exi_Proposition_}).

This is remarkable because a parallel statement
in K{\"a}hler geometry - whether any K{\"a}hler manifold
is a deformation of a projective one - is still
unknown. 

The manifold $\tilde M$ is not Stein; the existence 
of globally defined holomorphic functions on $\tilde M$
is a priori unclear. However, from the existence of 
$\Gamma'\subset \tilde G_\C$ with
$\tilde M / \Gamma'$  quasiregular,
we obtain that $\tilde M$ admits
quite a few holomorphic functions. 
These functions are obtained
{}from the sections of powers of
an ample bundle $L_Q$ on the leaf space
$M/\Xi$ (see Subsection \ref{_qr_Vaisman_Intro_Subsection_}).

Denote by $\gamma':\; \tilde M \arrow \tilde M$
the generator of $\Gamma'$, ang let $q\in \R^{>0}$ be the
monodromy of the weight bundle $L$. Consider the 
space $V$ of holomorphic functions
$\mu:\; \tilde M \arrow \C$, $\gamma'(\mu) = q^k\mu$.
This space is identified with the space of sections
of $L_Q^k$ (Subsection \ref{_qr_Vaisman_Intro_Subsection_}).
We shall think of $V$ as of eigenspace of $\gamma'$
acting $\calo_{\tilde M}$.

Clearly, $M$ is equipped with a natural map 
\begin{equation}\label{_M_to_quo_Equation_}
M \arrow V/\langle q^k \rangle.
\end{equation}
For $k$ sufficiently big, the functions
{}from $V$ have no common zeros, and the map
\eqref{_M_to_quo_Equation_} induces a holomorphic
immersion $M \arrow (V\backslash 0)/\langle q^k \rangle$
(\ref{_qr_imme_Intro_Theorem_}).

We obtain that the natural map $\tilde M \arrow (V\backslash 0)$
is an immersion. Denote by $\gamma$ the generator of $\Gamma$.
By construction, the maps $\gamma, \gamma': \tilde M \arrow \tilde M$
belong to the complexification $\tilde G_\C$. Since
$\tilde G$ is commutative, $\gamma$ and $\gamma'$ commute.
Since $V$ is an eigenspace of $\gamma'$, $\gamma$ 
acts on $V$. This gives a map
\[ M = \tilde M /\langle\gamma\rangle \arrow 
  (V\backslash 0)/\langle \gamma \rangle.
\]
This map is an immersion.

This way we obtain an immersion from arbitrary Vaisman 
manifold to a Hopf manifold $(V\backslash 0)/\Z$;
see \ref{_imme_the_gene_Theorem_} for details.

This Immersion Theorem is a Vaisman-geometric
analogue of the Kodaira embedding theorem.

%%%%%%%%%%%%%%%%%%%%%%%%%%%%%%%%%%%%%%%%%%%%%%%%
\subsection{Category of Vaisman varieties}
%%%%%%%%%%%%%%%%%%%%%%%%%%%%%%%%%%%%%%%%%%%%%%%%

Vaisman manifolds naturally form a category. In 
\cite{_Tsukada_} and \cite{_Verbitsky:LCHK_} it was shown that any closed
submanifold $X\subset M$ of a compact Vaisman 
manifold is again Vaisman, assuming $\dim X>0$.

Define a morphism of Vaisman manifolds as a 
holomorphic map commuting with the holomorphic Lee
flow. The Immersion Theorem can be stated as follows:

\hfill

%%%%%%%%%%%%%%%%%%%%%%%%%%%%%%%%%%%%
\theorem
Any Vaisman manifold admits a Vaisman immersion
to a Hopf manifold.

\hfill

\noindent
{\bf Proof:} See \ref{_imme_the_gene_Theorem_}.
\endproof

\hfill

We could define Vaisman varieties as varieties
immersed  in a Hopf manifold and equipped with the induced Lee
flow, and define morphisms of Vaisman varieties
as a holomorphic map commuting with the Lee flow.
This is similar to defining projective varieties 
as complex subvarieties in a complex projective space.

The Vaisman geometry becomes then a legitimate chapter
of complex algebraic geometry. In fact, the
data needed to define a Vaisman submanifold
in a Hopf manifold (a space $V$ with a linear
flow $\psi_t = e^{t A}, \  A \in {\goth {gl}}(V)$
and an analytic subvariety $X\subset V$) are
classical and well known. 

The Immersion Theorem suggests that Vaisman varieties
have interesting intrinsic geometry, in many ways
parallel to the usual algebraic geometry. One could ask,
for instance, the following.

\hfill

%%%%%%%%%%%%%%%%%%%%%%%%%%%%%%%%%%%%
\question
What are the moduli spaces of Vaisman manifolds?
They are, generally speaking, 
distinct from the moduli of complex 
analytic deformations: not every complex
analytic deformation is Vaisman.

\hfill

%%%%%%%%%%%%%%%%%%%%%%%%%%%%%%%%%%%%%%%%%%%%%%%%
\question
Is there a notion of stable holomorphic bundle
on a Vaisman manifold? What are the moduli of holomorphic
bundles? 

\hfill

%%%%%%%%%%%%%%%%%%%%%%%%%%%%%%%%%%%%%%%%%%%%%%%%
\question
Given a singular Vaisman variety (that is, a subvariety 
of a Hopf manifold), do we have a resolution of singularities
within Vaisman category? Given a holomorphic map 
of Vaisman varieties, is it always compatible 
with the Lee field? What are the birational maps 
of Vaisman varieties?

\hfill

%%%%%%%%%%%%%%%%%%%%%%%%%%%%%%%%%%%%%%%%%%%%%%%%
\question 
Consider a Vaisman manifold $M$, $\dim M =n$,
with canonical bundle isomorphic to $n$-th degree
of a weight bundle. Is there an 
Einstein-Weyl\footnote{Please consult 
\cite{_Calderbank_Pedersen_} or \cite{_Ornea:LCHK_}
for a reference to Eistein-Weyl LCK structures.}
LCK metric on $M$? This will amount to
a Vaisman analogue of the Calabi-Yau theorem.

%%%%%%%%%%%%%%%%%%%%%%%%%%%%%%%%%%%%%%%%%%%%%%%%

\section{LCK manifolds}

%%%%%%%%%%%%%%%%%%%%%%%%%%%%%%%%%%%%%%%%%%%%%%%%

In this Section we give a brief exposition of Vaisman geometry, 
following \cite{_Dragomir_Ornea_}, \cite{_Kamishima_Ornea_}, 
\cite{_Verbitsky:LCHK_}.

%%%%%%%%%%%%%%%%%%%%%%%%%%%%%%%%%%%%%%%%%%%%%%%%%%%%%%%%%%%%
\subsection{Differential geometry of LCK manifolds}
%%%%%%%%%%%%%%%%%%%%%%%%%%%%%%%%%%%%%%%%%%%%%%%%

Throughout this paper $(M,I,g)$ will denote a connected Hermitian
manifold of complex dimension $m\geq 2$, with fundamental two-form
defined by $\omega(X,Y)=g(X,IY)$.  

%%%%%%%%%%%%%%%%%%%%%%%%%%%%%%%%%%%%
\hfill

\definition\label{_LCK_Definition_}
$(M,I,g)$ is called {\bf
 locally conformal   K{\"a}hler},  {\bf LCK} for short,
if there exists a global {\em closed} one-form $\theta$ satisfying the
following integrability condition:
\begin{equation}\label{1}
d\omega=\theta\wedge  \omega.
\end{equation}

\hfill

The one-form $\theta$ is called the {\bf Lee form} and its metric dual
field, denoted $\theta ^\sharp$, is called the {\bf Lee field}. Locally we
may write $\theta=df$, hence the local metric $e ^{-f}g$ is K{\"a}hler,
thus motivating the definition.

\hfill

Te pull-back of $\theta$ becomes exact on the universal covering
space $\hat M$. The deck transformations of $\hat M$ induce
homotheties of the K{\"a}hler manifold $(\hat M, e ^{-f}g)$.
We obtain that $(M,I,g)$ is LCK
  if and only if its universal covering 
admits a K{\"a}hler metric, with deck transformations
acting by homotheties. 
This is in fact true for any covering
  $\tilde M$ such that the pullback of $\theta$ is exact
on $\tilde M$. 
For such a covering we  let $\mathcal{H}(\tilde M)$ denote the group
  of all holomorphic homotheties of $\tilde M$ and let 
\begin{equation}\label{_rho_}
\rho:\mathcal{H}(\tilde M)\arrow \R^{>0}
\end{equation}
be the homomorphism which associates to each homothety its scale
factor. 

\hfill

LCK manifolds are examples of what is called {\bf a Weyl manifold}.
Recall that Weyl manifold is a Riemannian 
manifold $(M, g)$ admitting a torsion-free
connection which satisfies 
\begin{equation}\label{_Weyl_Equation_} 
 \nabla(g) = g\otimes \theta,
\end{equation}
or, equivalently, preserves the conformal class of $g$.
For a detailed overview of Weyl geometry, see
\cite{_Calderbank_Pedersen_} and \cite{_Ornea:LCHK_}.
A Weyl connection on an LCK manifold is provided
by the Levi-Civita connection $\nabla_{LC}$ on its K{\"a}hler
covering $\tilde M$. The deck transforms of $\tilde M$
multiply the K{\"a}hler metric by constant, hence
commute with $\nabla_{LC}$. This allows one to pushdown
$\nabla_{LC}$ to $M$. Clearly, $\nabla_{LC}$ preserves the 
conformal class of $g$, hence satisfies \eqref{_Weyl_Equation_}.

\hfill

%%%%%%%%%%%%%%%%%%%%%%%%%%%%%%%%%%%%%%%%%%%%%%%%%%%%%%%%%%%%
\definition
Let $(M, \nabla, g)$ be a Weyl manifold.
The real line bundle assciated to the  bundle of linear frames of $(M,g)$ by 
the
representation $GL(n,\R)\ni A\mapsto |\det A|^{\frac 1n}$ is called
the {\bf weight bundle} of $M$.

\hfill

The Weyl connection
induces a flat connection, also denoted  $\nabla$, in $L$. Its
connection 1-form can be identified with $\theta$, hence $(L, \nabla)$
is a flat bundle.

Usually one has no
way to choose a specific metric in a conformal class, hence to
trivialize $L$. But on compact LCK manifolds one has such a
possibility due to a result of Gauduchon:

\hfill

%%%%%%%%%%%%%%%%%%%%%%%%%%%%%%%%%%%%%%%%%%%%%%%%
\claim \label{_Gaud_metric_Claim_}
\cite{_Gauduchon_1984_} 
Let $([g], \nabla)$ be a Weyl structure
on a  compact manifold.  There
exists a unique (up to homothety) metric in $[g]$ with coclosed Lee
form. 

\hfill

Clearly, any metric conformal with a LCK one is still LCK (with
respect to a fixed complex structure), so that we can give:   

\hfill

%%%%%%%%%%%%%%%%%%%%%%%%%%%%%%%%%%%%%%%%%%%%%%%%%%%%%%%%%%%%
\definition
 The unique (up to homothety) LCK metric in the conformal class $[g]$
 of $(M,I,g)$ with harmonic associated Lee
 form is called the {\bf Gauduchon metric}.

\hfill

Note that the Gauduchon metric endowes $L$ with a trivialization $l$. 

%%%%%%%%%%%%%%%%%%%%%%%%%%%%%%%%%%%%%%%%%%%%%%%%
\subsection{Vaisman manifolds}
\label{_Vaisman_Subsection_}
%%%%%%%%%%%%%%%%%%%%%%%%%%%%%%%%%%%%%%%%%%%%%%%%

The typical example of LCK manifold, is
the diagonal Hopf manifold $H_{\alpha}:=(\C^{n}\setminus 0)/\Z$ where $\Z$ is 
generated
by $z\mapsto \alpha z$, $|\alpha|\neq 0,1$. The LCK metric is here
$|z|^{-2}\sum dz_i\otimes d\bar z_i$ with Lee form given (locally) by $-d\log 
|z|^2$ and one can check that it is parallel with
respect to the Levi Civita connection of the LCK metric. Manifolds
with such structure were intensively studied by  Vaisman under the
name of generalized Hopf manifolds which later proved to be
inappropriate. So that we adopt the following:

\hfill 

%%%%%%%%%%%%%%%%%%%%%%%%%%%%%%%%%%%%%%%%%%%%%%%%%%%%%%%%%%%%
\definition \label{_Vaisman_Definition_}
A LCK manifold $(M,I,g)$ with parallel Lee form is called a
 {\bf Vaisman manifold}.

\hfill

Every Hopf surface admits a LCK metric, \cite{_Gauduchon_Ornea_}, but
the Hopf surfaces with K{\"a}hler rank 0 do not admit Vaisman metrics,
\cite{_Belgun_}. The construction in  \cite{_Gauduchon_Ornea_} was
generalized in \cite{_Kamishima_Ornea_} where Vaisman metrics were
found on manifolds $H_\Lambda:=(\C^n\setminus 0)/\Z$ with
$\Lambda = (\alpha_1,\ldots,\alpha_n)$, $\alpha_i\in \C$, $|\alpha_1|\geq 
\cdots\geq |\alpha_n|>1$ and $\Z$ 
generated by the transformations $(z_j)\mapsto (\alpha_jz_j)$. 

\hfill

%%%%%%%%%%%%%%%%%%%%%%%%%%%%%%%%%%%%%%%%%%%%%%%%%%%%%%%%%%%%
\definition \label{_Hopf_mfld_Definition_}
A linear operator on a complex vector space $V$ 
is called {\bf semisimple} if it can be diagonalized.
Given a semisimple operator $A:\; ÷ \arrow B$ 
with all eigenvalues  $|\alpha_i|>1$, the quotient
$H_\Lambda:= (V\backslash 0)\langle A \rangle$ 
described above is called {\bf a Hopf manifold}.

\hfill

The Vaisman structure on $H_\Lambda$ is intimately related by
a diffeomorphism between this manifold and $S^1\times S^{2n-1}$ to a 
specific  Sasakian structure of the odd sphere. For an up to date
introduction in Sasakian geometry we refer to \cite{_Boyer_Galicki_}~;
here we just recall that $(X,h)$ is Sasakian if and only if the cone
$(\mathcal{C}(X):=X\times \R^{>0}, dt^2+t^2h)$ is K{\"a}hlerian. The
example of the Hopf manifold is not casual: the total space of any
flat $S^1$ bundle over a compact Sasakian manifold can be given a
Vaisman structure. Recently, the following structure theorem was
proved:

\hfill
  
\theorem \cite{_OV:Structure_} Any compact Vaisman manifold is a Riemannian 
sus\-pen\-sion with Sa\-sa\-kian
fibers over a circle,  and conversely, the
total space of such a Riemannian suspension is a compact Vaisman
manifold. 

\hfill

The proof of this structure theorem relies on the following
particular property of the weight bundle that we shall use also in the
present paper:

\hfill

\theorem \cite{_OV:Structure_} The weight bundle of a compact Vaisman 
manifold has monodromy
$\Z$.

\hfill

An LCK metric with parallel Lee form on a compact manifold necessarily
coincides with a Gauduchon metric. Moreover, on compact Vaisman
manifolds $(M,I,g)$, the group $\Aut(M)$ of all conformal
biholomorphisms coincides with $\Iso(M,g)$,
\cite{mpps}, the isometry group of the Gauduchon metric, thus being
compact.

\hfill

We recall that on a Vaisman manifold, the flow of the Lee field (which
we call Lee flow) is formed by holomorphic isometries of the Gauduchon
metric. The same is true for the flow of the anti-Lee field $I\theta
^\sharp$.  We shall denote by $G$ the real Lie subgroup of $\Aut(M)$
generated by the Lee flow.

\hfill

\definition
The covering  $\pi:\; \tilde M\arrow M$   
associated with the monodromy group of $L$ is called the {\bf weight
covering} of $M$.

\hfill

We let $\Aut(\tilde M, M)$ be  
the group of all conformal automorphisms of $\tilde M$
which make the following diagram commutative: 

\begin{equation}
\begin{CD}
 \tilde M@>{\tilde f}>> \tilde M \\
@V{\pi}VV  @VV{\pi}V              \\
M@>{f}>>  M 
\end{CD}
\end{equation}

By the above, $\Aut(\tilde M, M)$ is nonempty, as it contains the Lee flow. 
For further use let 
\begin{equation}\label{_M_automo_forgetful_Equation_} 
 \Phi:\; \Aut (\tilde M, M)\arrow \Aut (M)
\end{equation}
be the forgetful map.

\hfill

Much of the geometry of a LCK manifold, especially of a Vaisman
manifold, can be expressed in terms of the two-form
$\omega_0=d^c\theta$, where $d^c=-IdI$, introduced and studied in
\cite{_Verbitsky:LCHK_}. It is semipositive, its only zero eigenvalue
being in the direction of the Lee field. To further understand
$\omega_0$, consider the complexified bundle $L_\C:= L\otimes_\R \C$. The 
$(0,1)$ part
of the complexified Weyl connection endowes this complex
line bundle  with a holomorphic structure. We also equip $L_\C$ with a
Hermitian structure such that to normalize to 1 the length of the 
trivialization
$l$. Whenever we refer to $L_\C$, we implicitely refer to this
holomorphic and Hermitian structure. In this setting, the form $\omega_0$ has 
the following geometric meaning:

\hfill

\claim \cite[Th. 6.7]{_Verbitsky:LCHK_}\label{_curva_L_omega_0_Claim_}
The curvature of the Chern connection of $L_\C$ is equal to  $ -2 \1\omega_0$.

%Let $M$ be an LCK manifold, and $L$ the 
%corresponding weight bundle equipped with a 
%holomorphic and a Hermitian structure as above. Denote by
%${}^C\nabla$ the standard Hermitian connection on $L$
%(so-called Chern connection), and let $C$ be its
%curvature. Then $C= -2 \1\omega_0$, where 
%$\omega_0 = d^c \theta$ is the 
%standard 2-form on $M$.

%{\bf Proof:} This is \cite{_Verbitsky:LCHK_}, Theorem 6.7.
%\endproof

\hfill

We end this section by recalling from \cite[Pr. 4.4]{_Verbitsky:LCHK_} that 
if the Lee form of a Vaisman 
manifold is exact, $\theta=df$, then the function $r:=e ^{-f}$ is a
potential for the K{\"a}hler metric $r\omega$, \emph{i.e.} 
\begin{equation}\label{_kahler_pote_Equation_} 
r\omega=dd^cr.
\end{equation}

%%%%%%%%%%%%%%%%%%%%%%%%%%%%%%%%%%%%%%%%%%%%%%%%%%%%%%%%%%%%

\section{Quasiregular  Vaisman manifolds}

%%%%%%%%%%%%%%%%%%%%%%%%%%%%%%%%%%%%%%%%%%%%%%%%%%%%%%%%%%%%

%%%%%%%%%%%%%%%%%%%%%%%%%%%%%%%%%%%%%%%%%%%%%%%%%%%%%%%%%%%%
\subsection{Quasiregular  Vaisman manifolds and the weight bundle}
\label{_qr_Vais_Subsection_}
%%%%%%%%%%%%%%%%%%%%%%%%%%%%%%%%%%%%%%%%%%%%%%%%%%%%%%%%%%%%

As the Lee field is Killing and holomorphic on a Vaisman manifold, the
distribution locally generated by $\theta ^\sharp$ and $I\theta
^\sharp$ defines a holomorphic, Riemannian, totally geodesic foliation
$\Xi$ that we call {\bf the  Lee foliation}. Its leaves are elliptic
curves (cf. \cite{_Vaisman:Dedicata_}).\footnote{Some algebraic
geometers prefer to distinguish between ``compact 1-dimensional
complex tori'' and ``elliptic curves'' (compact 1-dimensional 
complex tori with a marked point). The compact fibers of the 
Lee foliation are, in this terminology, compact 1-dimensional
complex tori.} The leaves of $\Xi$ are identified
with the orbits of a holomorphic flow generated by
$\theta ^\sharp$; this flow is called {\bf the complex
Lee flow}.

\hfill

\definition
A Vaisman manifold is called:
\begin{itemize}
\item  {\bf quasiregular } if $\Xi$ has compact leaves. In this case
  the leaf space $M/\Xi$ is a Hausdorff orbifold. 
\item  {\bf regular} if it is quasiregular and 
the quotient map $M \arrow M/\Xi$ is smooth.
\end{itemize}

\hfill

Combining results from \cite{_Vaisman:Dedicata_} and
\cite{_Verbitsky:LCHK_}, we can state:

\hfill

%%%%%%%%%%%%%%%%%%%%%%%%%%%%%%%%%%%%%%%%%%%%%%%%
\claim\label{_Lee_fibra_Claim_}
Let $M$ be a quasiregular  Vaisman manifold 
and let $Q$ be the leaf space of the Lee foliation.
Then $Q$ is a K{\"a}hler orbifold and the pull-back of its K{\"a}hler form
by the natural projection coincides with $\omega_0$.

\hfill
  
%%%%%%%%%%%%%%%%%%%%%%%%%%%%%%%%%%%%%%%%%%%%%%%%
\definition
Let $M$ be a quasiregular  Vaisman manifold and denote  
$f: M \arrow Q$ the above associated fibration. Then $f: M \arrow Q$
is called {\bf the Lee fibration of $M$}.

\hfill

Let $f_*L_\C\arrow Q$ be the push-forward of the bundle $L_\C$. The main 
result of this Subsection is the following Theorem.

\hfill

%%%%%%%%%%%%%%%%%%%%%%%%%%%%%%%%%%%%%%%%%%%%%%%%%%%%%%%%%%%%
\theorem \label{_L_trivi_on_fibe_Theorem_}
Let $M$ be a compact quasiregular  Vaisman manifold,
and let $f:\; M \arrow Q$ its Lee fibration.
% Consider theweight bundle $L$ as a holomorphic line bundle on $M$. 
Then
$L_\C$ is trivial along the fibers of $f$: the natural map
\begin{equation}\label{_triviality_on_fibe_Equation_}
L_\C \arrow f^* f_* L_\C 
\end{equation}
is an isomorphism.

\noindent{\bf Proof:} The statement
of \ref{_L_trivi_on_fibe_Theorem_} is local on $Q$.
Since $f$ is proper, it suffices to prove that
$L_\C$ is trivial on $f^{-1}(U)$, for a sufficiently
small $U\subset Q$. Indeed, for any 
proper map with connected fibers, a 
pushforward of a structure sheaf is a structure
sheaf, hence the equation 
\eqref{_triviality_on_fibe_Equation_}
is true for trivial bundles. 

Let $\tilde M$ be the weight covering of $M$. The lifting $\tilde L_\C$
of $L_\C$ to $\tilde M$ is flat and has trivial monodromy.  
Therefore, $\tilde L_\C$ is trivialized, in a canonical way. 
This allows us to interpret sections of $\tilde L_\C$ as 
holomorphic functions on $\tilde M$. A function $\mu:\; \tilde M \arrow \C$
of $\tilde L_\C$ is lifted from $L_\C$ if and only
if it satisfies the following monodromy condition:
\begin{equation}\label{_monodro_on_secti_L_Equation_}
\gamma(\mu) = q\mu,
\end{equation}
where $\gamma:\; \tilde M \arrow \tilde M$ is the
generator of the deck transformation group of $\tilde M$, and
$q\in \End(L)=\R^{>0}$ the corresponding monodromy
action.

The fibers of $f$ are elliptic curves and coincide with the orbits of
$G_\C$, the group generated by the complex Lee flow. 
For any orbit $C\subset M$, the corresponding covering
$\tilde C \subset \tilde M$ is isomorphic to $\C^*$,
with $\gamma$ acting on $\tilde C \cong \C^*$ as a multiplication
by $q$. After this identification,
the equation \eqref{_monodro_on_secti_L_Equation_}
is rewritten as 
\begin{equation}\label{_mono_on_secti_on_fibers_Equation_}
\mu(qz) = q \mu(z), \ \ z\in \tilde C =\C^*.
\end{equation}
A function on $\C^*$ satisfying \eqref{_mono_on_secti_on_fibers_Equation_}
must be linear. We obtain that the Lee field acts linearly on
all functions satisfying \eqref{_monodro_on_secti_L_Equation_}:
\begin{equation}\label{_monodro_ODE_Equation_}
\Lie_{\theta^\sharp}(\mu) = \log q \cdot \mu.
\end{equation}
This is an ordinary differential equation.
Using the uniqueness and existence of solutions of ODE, we
arrive at the following claim.

\hfill

%%%%%%%%%%%%%%%%%%%%%%%%%%%%%%%%%%%%%%%%%%%%%%%%
\claim\label{_L_trivia_with_slice_Claim_}
Let $M$ be a quasiregular Vaisman manifold,
$G_\C$ the complex Lie group generated by the Lie
flow, $f: M \arrow Q$ the quotient space,
$Q=M/G_\C$, and $U\subset M$ the {\bf slice}
of $f$, that is, a submanifold of $M$
(not necessarily closed in $M$) such that the
restriction $f\restrict U:\; U \arrow Q$ is
an open embedding. Consider the Vaisman
manifold $M_U:= f^{-1}(f(U))$.
Pick a trivialization $\mu_U$
of $L_\C\restrict U$ (this is possible
to do if $U$ is sufficiently small).
This gives a trivialization of $L_\C$ on $M_U$.

\hfill

\noindent{\bf Proof:}
We solve the ODE \eqref{_monodro_ODE_Equation_}
with an initial  condition $\mu\restrict U = \mu_U$.
The solution gives a trivialization of $L_\C$
as the above argument indicates.
\endproof

\hfill

Return to the proof of \ref{_L_trivi_on_fibe_Theorem_}.
A slice of the action of $G_\C$ exists at any point $x\in M$ 
where $G_\C$ acts smoothly, that is, for all points $x\in M$ 
which are not critical points of $f:\; M \arrow Q$. 
Using \ref{_L_trivia_with_slice_Claim_}, we obtain
a trivialization needed in \ref{_L_trivi_on_fibe_Theorem_}
for all regular Vaisman manifolds (regular Vaisman
manifolds are precisely those manifolds for which
$f$ is smooth). When $M$ is only quasiregular,
we use the orbifold covering to reduce \ref{_L_trivi_on_fibe_Theorem_}
to regular case, as follows.

Let $x\in Q$ be a critical value of $f$.
Since \ref{_L_trivi_on_fibe_Theorem_} is local,
we may replace $Q$ with a sufficiently small
neighbourhood of $U\ni x$, and $M$ with $f^{-1}(U)$.
Since $f$ is an orbifold morphism, it admits a finite
covering by a smooth morphism of manifolds $M'\stackrel{f'} \arrow Q'$
such that the horizontal maps of the Cartesian square
\[
\begin{CD}
M' @>{\tau_M}>> M\\
@V{f'}VV  @VV{f}V \\
Q'@>{\tau}>>Q\\
\end{CD}
\]
are etale in the orbifold category\footnote{Locally we have
  $\tau:V'/G'\arrow V/G$ for some open sets $V$, $V'$ and finite
  groups $G'$, $G$. Moreover, $V$ can be chosen simply connected. By
  etale in the orbifold category we understand that the induced map
  $V\arrow V'$ is etale}.
The Vaisman manifold $M'$ is by construction regular.
Choosing $Q$ sufficiently small, we may insure that
$f'$ admits  a slice satisfying the assumptions of 
\ref{_L_trivia_with_slice_Claim_}.
Then the lift $L'_\C$ of $L_\C$ to 
$M'$ is a trivial bundle. Denote by $\cal G$ the
Galois group of $\tau$, that is, the orbifold
deck transform group of the covering
$\tau:\; Q' \arrow Q$. The group $\cal G$
acts on the pushforward ${\tau_M}_* L'_\C$,
and we have
\[ L_\C = ({\tau_M}_* L'_\C)^{\cal G},\]
where $(\cdot)^{\cal G}$ denotes the sheaf
of $\cal G$-invariants. Since $L'_\C$ is trivial,
the bundle $L_\C$ is isomorphic to 
$({\tau_M}_* \calo_{M'})^{\cal G}$.
This bundle is trivial, because 
$\tau_M$ is etale. We proved that
 $L_\C$ is trivial, locally in $Q$.
\ref{_L_trivi_on_fibe_Theorem_} is proven.
The same argument shows that the natural Hermitian structure
on $L_\C \cong f^* f_*L_\C$ is lifted from the Hermitian structure
on the bundle $f_*L_\C$.
\endproof

%%%%%%%%%%%%%%%%%%%%%%%%%%%%%%%%%%%%%%%%%%%%%%%%
\subsection[Kodaira-Nakano theorem for quasiregular Vaisman 
manifolds]{Kodaira-Nakano theorem for quasiregular \\  Vaisman manifolds}
%%%%%%%%%%%%%%%%%%%%%%%%%%%%%%%%%%%%%%%%%%%%%%%%

Let $M$ be a compact quasiregular  Vaisman manifold, $L_\C$ the weight bundle,
and $C = -2 \1\omega_0$ the curvature of the Chern connection
in $L_\C$ (\ref{_curva_L_omega_0_Claim_}). It was shown in
\cite{_Verbitsky:LCHK_} that $\omega_0$ is obtained as a 
pullback of a K{\"a}hler form on $Q$:
\[ \omega_0 = f^* \omega_Q.
\]

\hfill

%%%%%%%%%%%%%%%%%%%%%%%%%%%%%%%%%%%%%%%%%%%%%%%%
\proposition\label{_L_Q_ample_Proposition_}
Let $M$ be a quasiregular  Vaisman manifold,
$f:\; M \arrow Q$ the Lee fibration, $L_\C$ the weight bundle
and $L_Q:= f_*L_\C$ the corresponding line bundle on $Q$, equipped
with the Chern connection.  Then 
\[ 
C_Q= -2 \1\omega_Q.
\]
In particular, $L_Q$ is ample.

\hfill

\noindent{\bf Proof:} The Hermitian bundle $L_Q$ is equipped
with a natural Hermitian metric $h$, such that
 the metric on $L_\C=f^*L_Q$ is a pullback of $h$
(see \ref{_L_trivi_on_fibe_Theorem_}). 
Therefore, the Chern connection
$\nabla_C$ on $L_\C$ is lifted from $L_Q$. The curvature of $\nabla_C$
is $-2 \1\omega_0$, as \ref{_curva_L_omega_0_Claim_} implies. 
The form $-2 \1\omega_0$ is a pullback of $-2 \1\omega_Q$.
(\ref{_Lee_fibra_Claim_}). This proves \ref{_L_Q_ample_Proposition_}.
\endproof

\hfill

{}From \ref{_L_Q_ample_Proposition_} and the 
Kodaira embedding theorem for orbifolds (cf. \cite{_Baily_}), we obtain the
following result

\hfill  

%%%%%%%%%%%%%%%%%%%%%%%%%%%%%%%%%%%%%%%%%%%%%%%%
\corollary \label{_KN_for_qr_Vaisman_Corollary_}
Let $M$ be a compact quasiregular  Vaisman manifold,
$f:\; M \arrow Q$ the Lee fibration and $L_\C$ the weight bundle.
Then, for $k$ sufficiently big, 
the bundle $L_\C^k$ has no base points, and defines a natural map
$l_k :\; M \arrow {\C P}^{n-1}$, in a usual way
($\dim H^0(L_\C^k) = n$). Moreover, 
$l_k$ is factorized through $f:\; M \arrow Q$:
\[  l_k = \underline{l}_k\circ f, 
    \ \ \underline{l}_k:\; Q \arrow {\C P}^{n-1}.
\]
and $\underline{l}_k:\; Q \arrow {\C P}^{n-1}$ is 
an embedding.

%%%%%%%%%%%%%%%%%%%%%%%%%%%%%%%%%%%%%%%%%%%%%%%%
\subsection{Immersion theorem for quasiregular  Vaisman manifolds}
\label{_qr_imme_Subsection_}
%%%%%%%%%%%%%%%%%%%%%%%%%%%%%%%%%%%%%%%%%%%%%%%%

The Kodaira theorem (\ref{_KN_for_qr_Vaisman_Corollary_})
implies the following Immersion Theorem for quasiregular  Vaisman
manifolds.

Recall that the weight bundle $L$ is equipped with 
a flat connection. Consider the weight covering
$\tilde M \arrow M$ associated with the 
monodromy group of $L$. By definition,
$\tilde M$ is the smallest covering 
such that the pullback $\tilde L$ of
$L$ to $\tilde M$ has trivial monodromy.
Since $\tilde L$ is a flat bundle with
trivial monodromy, it has a canonical
trivialization.

Let $l_1$, $l_2$, ... ,
$l_n$ be a basis in the space $H^0(L_\C^k)$.
Since $\tilde L$ is equipped with a  canonical
trivialization, $l_i$ gives a function
$\lambda_i:\; \tilde M \arrow \C$.
Since $l_1$, $l_2$, ... ,
$l_n$ induce an embedding $Q \hookrightarrow \C P^{n-1}$,
the functions $\lambda_i$ have no common zeros.
Together these function define a map
\begin{equation}\label{_tilde_lambda_defi_Equation_} 
  \tilde \lambda:\; \tilde M \arrow \C^n \backslash 0.
\end{equation}
The monodromy group $\Gamma$ of $L$ is isomorphic to $\Bbb Z$
(\cite{_OV:Structure_}; see also Subsection
\ref{_tilde_G_intro_Subsection_}).  Denote by $\gamma$ its generator,
and let $q\in \Aut(L)\cong \R^{>0}$ be the action of $\gamma$
on the fibers of $L_\C$. A holomorphic section $\tilde l$ of $\tilde L_\C^k$
is obtained as a pullback from a section of $L_\C^k$ if 
and only if
\begin{equation}\label{_gamma_acts_on_pullbacks_Equation_}
\gamma(\tilde l) = q^k \tilde l,
\end{equation}
where $\gamma:\; H^0(\tilde L_\C^k) \arrow H^0(\tilde L_\C^k)$
is the natural equivariant action of $\gamma\in \Gamma$ on
$\tilde L_\C^k$.

Now return to the map \eqref{_tilde_lambda_defi_Equation_}.
The eigenspace property \eqref{_gamma_acts_on_pullbacks_Equation_}
implies that the following square is commutative
\begin{equation}\label{_tilde_M_immersion_CD_Equation_}
\begin{CD}
\tilde M @>{\tilde \lambda}>> \C^n \backslash 0\\
@V{\gamma}VV  @VV{\text{mult. by $q^k$}}V \\
\tilde M @>{\tilde \lambda}>>\C^n \backslash 0\\
\end{CD}
\end{equation}
Therefore, $\tilde \lambda$ is obtained as a covering for a map
$\lambda:\; M \arrow \left(\C^n \backslash 0 \big / \langle 
q^k\rangle\right)$ to
a Hopf manifold. We also have a commutative square
of Lee fibrations
\begin{equation}\label{_Lee_fibra_commu_Equation_}
\begin{CD}
M @>{\lambda}>> \left(\C^n \backslash 0 \big / \langle q^k\rangle\right)\\
@V{f}VV  @VVV \\
Q @>{\underline l_k}>> \C P^{n-1}.
\end{CD}
\end{equation}
The map $\underline l_k$ is by construction
an embedding, hence $\lambda$ maps different fibers of $f$
to different Lee fibers of the Hopf manifold 
$\left(\C^n \backslash 0 \big / \langle q^k\rangle\right)$.
Both sides of \eqref{_Lee_fibra_commu_Equation_} are 
Lee fibrations, with fibers isomorphic to 1-dimensional 
compact tori, with the affine coordinates provided
by the sections of the weight bundle.
The map $\lambda$ is by construction compatible with 
affine structure on these fibers.
Therefore, $\lambda$ induces a finite
covering on the fibers of Lee fibrations of 
\eqref{_Lee_fibra_commu_Equation_}.
We obtain the following theorem.

\hfill

%%%%%%%%%%%%%%%%%%%%%%%%%%%%%%%%%%%%%%%%%%%%%%%%
\theorem\label{_imme_for_qr_Vaisman_Theorem_}
Let $M$ be a compact quasiregular  Vaisman bundle. Then 
$M$ admits a holomorphic immersion
$\lambda:\; M \arrow \left(\C^n \backslash 0 \big / \langle
q^k\rangle\right)$ to
a (diagonal) Hopf manifold. Moreover, $\lambda$ is etale to its image,
that is, induces a finite covering from $M$ to $\lambda(M)$.

\hfill

In the rest of this paper we generalize this result to 
arbitrary (not necessarily quasiregular ) compact Vaisman
manifolds.

%%%%%%%%%%%%%%%%%%%%%%%%%%%%%%%%%%%%%%%%%%%%%%%%%%%%%%%%%%%%

\section{The group generated by the Lee flow and its applications}

%%%%%%%%%%%%%%%%%%%%%%%%%%%%%%%%%%%%%%%%%%%%%%%%%%%%%%%%%%%%

%%%%%%%%%%%%%%%%%%%%%%%%%%%%%%%%%%%%%%%%%%%%%%%%
\subsection{The groups $G$, $\tilde G$ generated by the Lee flow}
\label{_tilde_G_intro_Subsection_}
%%%%%%%%%%%%%%%%%%%%%%%%%%%%%%%%%%%%%%%%%%%%%%%%

%%%%%%%%%%%%%%%%%%%%%%%%%%%%%%%%%%%%%%%%%%%%%%%%
\definition
Let $G$ be the smallest Lie subgroup of $\Aut(M)$ containing the Lee
flow and let  $\tilde G$ be its lift to the weight covering $\tilde M$.

\hfill

It is easy to see that $G$ is compact and abelian
(see \cite{_Kamishima_Ornea_}). Therefore,
$G\cong (S^1)^k$.
In \cite{_OV:Structure_} it was shown that 
\begin{equation}\label{_tilde_G_structure_Equation_}
\tilde G \cong (S^1)^{k-1}\times \R,
\end{equation}
and the (restriction of) natural forgetful map $\tilde G \stackrel \Phi 
\arrow G$
 (cf. \eqref{_M_automo_forgetful_Equation_}) is a covering . The group 
$\Gamma$ of deck transformations of
$\tilde M$ is clearly isomorphic to $\ker \Phi$.
Therefore, \eqref{_tilde_G_structure_Equation_}
implies 
\begin{equation}
\Gamma \cong \Z.
\end{equation}

%%%%%%%%%%%%%%%%%%%%%%%%%%%%%%%%%%%%%%%%%%%%%%%%
\subsection{Complexification of $\tilde G$}
\label{_tilde_G_C_Subsection_}
%%%%%%%%%%%%%%%%%%%%%%%%%%%%%%%%%%%%%%%%%%%%%%%%

Consider the Lie algebra $\Lie(\tilde G)\subset H^0(T\tilde M)$.
The vector fields in $\Lie(\tilde G)$ are holomorphic,
because the Lee flow is holomorphic. 
Let $\tilde G_\C\subset \Aut(\tilde M)$
be the complex Lie group generated by $\tilde G$ and acting
on $\tilde M$. Clearly, 
$\Lie(\tilde G_\C)=\Lie(\tilde G)+I(\Lie(\tilde G))$,
where $I\in \End(T\tilde M)$ is the complex structure
operator on $\tilde M$.

\hfill

This operation is called {\bf complexification of
the Lie group $\tilde G$}. In the same way,
a complexification can be performed for any
real Lie group acting on a complex manifold
by holomorphic diffeomorphisms.

\hfill

By Remmert-Morimoto theorem,
a complex connected abelian Lie group
is isomorphic to
\begin{equation}\label{_Remmert_Morimoto_Equation_} 
(\C^*)^{l_1}\times \C^{l_2}\times T,
\end{equation}
where $T$ is a compact complex torus (\cite{_Morimoto_}). A Lie group
$\cal G$ is called {\bf linear } if  it has an exact complex
representation $\mathcal{ G}\hookrightarrow GL(n, \C)$.
Clearly, an abelian group $\cal G$ is linear  if 
$\mathcal{ G} \cong (\C^*)^{l_1}\times \C^{l_2}$.

\hfill

%%%%%%%%%%%%%%%%%%%%%%%%%%%%%%%%%%%%%%%%%%%%%%%%
\lemma
Let $M$ be a compact Vaisman manifold,
$\tilde M$ its weight covering, and 
$\tilde G_\C$ the complex Lie group 
generated by the Lee flow as above. Then 
$\tilde G_\C$ is linear .

\hfill

\noindent{\bf Proof:} Consider the action of $\tilde G_\C$
in $\tilde M$. Since $\tilde G_\C\subset \Diff(\tilde M)$, 
a general orbit $X$ of $\tilde G_\C$ is isomorphic to $\tilde G_\C$.
Since $\tilde M$ is equipped with a K{\"a}hler potential,
the manifold $X\cong \tilde G_\C$ is also equipped with
a K{\"a}hler potential. Therefore, $\tilde G_\C$ does
not contain a non-trivial compact complex torus.
\endproof

\hfill

%%%%%%%%%%%%%%%%%%%%%%%%%%%%%%%%%%%%%%%%%%%%%%%%%%%%%%%%%%%%
\proposition\label{_tilde_G_complexifi_Proposition_}
Let $M$ be a compact Vaisman manifold,
$\tilde M$ its weight covering, and 
$\tilde G_\C$ the complex Lie group 
generated by the Lee flow as above. Then 
$\tilde G_\C \cong (\C^*)^l$,
for some $l$. 

\hfill

\noindent{\bf Proof:} Denote by $\tilde G_0$ the group of all
$g \in \tilde G$ satisfying $\rho(g)=1$, where
$\rho:\; \tilde G \arrow \R^{>0}$ is the map defined in
\eqref{_rho_}. Since $\tilde G \cong (S^1)^{k-1}\times \R$,
and $\rho:\; \tilde G \arrow \R$ is surjective\footnote{The Lee flow acts 
by non-trivial homotheties on the K{\"a}hler form.}, 
we have $\tilde G_0 \cong (S^1)^{k-1}$. Denote by 
$\tilde G_K$ the Lie subgroup of $\tilde G_\C$ 
generated by $\tilde G_0$ and $e^{t I(\theta^\sharp)}$, 
$t\in \R$, where $I(\theta^\sharp)$ is the complex conjugate
of the Lee field. It is easy to see (see e.g. 
\cite[Pr. 4.3]{_Dragomir_Ornea_}) that $e^{t I(\theta^\sharp)}$
acts on $M$ by holomorphic isometries
and preserves the K{\"a}hler metric on $\tilde M$.\footnote{The
last statement is clear because $I\theta ^\sharp$ is Killing and if 
$\theta=df$, then $\mathcal{L}_{I\theta ^\sharp}(e ^{-f}g)=I\theta ^\sharp(e 
^{-f})g=-e ^{-f}\theta(I\theta ^\sharp)g=0$.}
Therefore, the forgetful map $\Phi:\; \tilde G_K \arrow \Aut(M)$
sends $\tilde G_K$ to isometries of $M$.

The following elementary claim is used to prove 
\ref{_tilde_G_complexifi_Proposition_}.

\hfill

%%%%%%%%%%%%%%%%%%%%%%%%%%%%%%%%%%%%%%%%%%%%%%%%
\claim\label{_tildeG_K_forgetful_Claim_}
Let $M$ be a compact Vaisman manifold,
\[ \tilde G_K\subset \Aut(\tilde M, M)\] the group
defined above, and 
\begin{equation}
\Phi:\; \tilde G_K \arrow \Iso(M)
\end{equation}
the forgetful map from $\tilde G_K$ to the group of 
isometries of $M$. Then $\Phi:\; \tilde G_K \arrow \Iso(M)$ 
is a monomorphism.

\hfill

\noindent{\bf Proof:}
Let $\nu \in \ker \Phi$. Since $\nu$ does not act on $M$,
$\nu$ belongs to the deck transformation group of $\tilde M$.
Unless $\nu$ is trivial, $\nu$ acts on $\tilde M$
by non-trivial homotheties of the K{\"a}hler form.
However, $\tilde G_K$ is generated by 
$\tilde G_0$ and $e^{t I(\theta^\sharp)}$,
and these diffeomorphisms preserve
the K{\"a}hler metric on $\tilde M$.
Therefore, the group $\tilde G_K$ preserves the
K{\"a}hler metric on $\tilde M$, and $\nu$
is trivial. \endproof

\hfill 

Return to the proof of \ref{_tilde_G_complexifi_Proposition_}.
The group $\tilde G\cong (S^1)^k \times \R$ is generated by 
$\tilde G_0\cong (S^1)^k$
and $e^{t \theta^\sharp}$, and $\tilde G_K$ is 
generated by $\tilde G_0$ and $e^{t I(\theta^\sharp)}$
Therefore the complexification of $\tilde G_K$ 
coincides with the complexification of $\tilde G$.

By \ref{_tildeG_K_forgetful_Claim_}, $\tilde G_K$
is a subgroup of the group of isometries
of $M$, and therefore it is compact.
A compact subgroup of a complex linear  Lie group
is {\bf totally real}, that is,
its Lie algebra does not intersect its complex
conjugate. Now, $\tilde G_\C$ is a
complexification of its totally real subgroup
 $\tilde G_K$; therefore, $\tilde G_K$
is a real form of $\tilde G_\C$.
A complex linear  connected abelian Lie group with 
compact real form is isomorphic
to $(\C^*)^n$, as the Morimoto-Remmert  classification
theorem implies (see \eqref{_Remmert_Morimoto_Equation_}). 
This proves \ref{_tilde_G_complexifi_Proposition_}.
\endproof

%%%%%%%%%%%%%%%%%%%%%%%%%%%%%%%%%%%%%%%%%%%%%%%%%%%%%%%%%%%%
\subsection{LCK manifolds associated with $\Gamma'\subset \tilde G_\C$}
%%%%%%%%%%%%%%%%%%%%%%%%%%%%%%%%%%%%%%%%%%%%%%%%%%%%%%%%%%%%

Let $M$ be a compact Vaisman manifold,
$\tilde M$ its weight covering, $G$, $\tilde G$, 
$\tilde G_\C$ the Lie groups defined in 
Subsection \ref{_tilde_G_intro_Subsection_}
and Subsection \ref{_tilde_G_C_Subsection_}.
By construction, $M = \tilde M/\Gamma$, where
$\Gamma\cong \Z$ is a subgroup of $\tilde G \subset \Aut(\tilde M,
M)$. Let $\gamma\in \Gamma$ be a generator of $\Gamma$.

\hfill

%%%%%%%%%%%%%%%%%%%%%%%%%%%%%%%%%%%%%%%%%%%%%%%%
\theorem\label{_defo_LCK_Theorem_}
In the above assumptions, consider 
an element $\gamma'\subset \tilde G_\C$ 
generating an abelian group 
$\Gamma'\subset \tilde G_\C$. Then,
for $\gamma'$ sufficiently close to $\gamma\in \tilde G_\C$,
the quotient space $\tilde M/\Gamma'$ is a compact
Vaisman manifold. Moreover, the Lee field of
$\tilde M/\Gamma'$ is proportional to
$v':=\log \gamma'\in \Lie(\tilde G_\C)$. 
\footnote{The group $\tilde G_\C$ is
connected and its exponential map is surjective, hence a definition of
the logarithm is possible.}

\hfill

\noindent{\bf Proof:} Let $v$ denote the Lee field $\theta^\sharp$ of
$M$ lifted to $\tilde M$.
The K{\"a}hler form on $\tilde M$ has a potential $\phi$
(see [Ve, Pr. 4.4] and the end of Section 2), which can be written as
$\phi := |v|^2$. If $\gamma$ is sufficiently
close to $\gamma'$, then $\phi':= |v'|^2$
is also plurisubharmonic and gives a K{\"a}hler 
potential on $\tilde M$. Denote by $\omega'$ the corresponding
K{\"a}hler  metric on $\tilde M$:
\[ \omega':= \1 \6\bar\6 |v'|^2
\]
Consider the differential flow $\psi'_t$ associated
with the vector field $v'$. Clearly, $\psi'_t= e^{tv'}$ multiplies
$v'$ by a number $e^t$. Therefore, $\psi'_t$ maps $|v'|^2$ 
to $e^{2t} |v'|^2$. We obtain that $\psi'_t$ preserves
the conformal class of the K{\"a}hler  metric $\omega'$. 
%\[ \omega'= \1 \6\bar\6 |v'|^2.
%\]
This implies that the quotient 
$\tilde M/\Gamma' = \tilde M/ \psi'_t$
is LCK. By construction, $v'$ is the Lee field
of $\tilde M/\Gamma'$, hence $\tilde M/\Gamma'$
is Vaisman. This proves \ref{_defo_LCK_Theorem_}.

We also obtain that $\tilde M/\Gamma'$ is Vaisman (in particular  LCK) 
whenever the function 
\[ 
\psi' := |\log \gamma'|^2
\]
is strictly plurisubharmonic.\footnote{A function $f$ is called
{\bf strictly plurisubharmonic} if and only if the 2-form
$\1 \6\bar\6 f$ is K{\"a}hler .}
\endproof

%%%%%%%%%%%%%%%%%%%%%%%%%%%%%%%%%%%%%%%%%%%%%%%%%%%%%%%%%%%%
\subsection{Quasiregular  Vaisman manifolds obtained as 
deformations}
\label{_qr_Vaisman_defo_Subsection_}
%%%%%%%%%%%%%%%%%%%%%%%%%%%%%%%%%%%%%%%%%%%%%%%%%%%%%%%%%%%%

Let $M$ be a compact Vaisman manifold
and $\tilde G_\C \cong (\C^*)^l$
the complex Lie group constructed in
Subsection \ref{_tilde_G_C_Subsection_}.
As we have shown, $\tilde G_\C$ has a compact
real form $\tilde G_K$. This allows as to equip
the Lie algebra 
\[ \Lie (\tilde G_\C)= \Lie (\tilde G_K)\otimes \C
\] 
with a rational structure, as follows.
Given $\delta\in \Lie (\tilde G_K),$
we say that $\delta$ is {\bf rational }
if the corresponding 1-parametric subgroup
$e^{t\delta}\subset \tilde G_K$
is compact. Clearly, this gives a rational
lattice in $\Lie (\tilde G_K)$,
and hence in $\Lie (\tilde G_\C)$.
This rational lattice corresponds
to an integer lattice of all 
\[ \{\delta\in \Lie (\tilde G_K) \ \  | \ \  \ e^\delta =1\}.
\]

A complex 1-parametric subgroup
$e^{\lambda \delta}\subset \tilde G_\C$
is isomorphic to $\C^*$ if and only if 
the line $\lambda \delta\subset \Lie (\tilde G_\C)$ 
is rational (contains a rational point); 
otherwise, $e^{\lambda \delta}$ is isomorphic to 
$\C$. This is clear, because the kernel of the 
exponential map consists exactly of integer elements
within $\Lie (\tilde G_K)$.

\hfill

%%%%%%%%%%%%%%%%%%%%%%%%%%%%%%%%%%%%%%%%%%%%%%%%%%%%%%%%%%%%
\proposition\label{_qr_defo_exi_Proposition_} 
Let $M$ be a compact Vaisman manifold, $\tilde M$
its weight covering and $G$, $\tilde G$, $\tilde G_\C$ the 
Lie groups defined above. Denote by $\Gamma\subset \tilde G$ the weight
monodromy group generated by $\gamma\in \tilde G$.
Take $\gamma'\in\tilde G_\C$, sufficiently
close to $\gamma$, in such a way that
the quotient $\tilde M/\langle \gamma'\rangle$ 
is a Vaisman manifold (see \ref{_defo_LCK_Theorem_}). Assume that
the line $\C \log \gamma'$ is rational.\footnote{Such $\gamma'$
are clearly dense in $\tilde G_\C$, because the set of rational
points is dense.} Then $\tilde M/\langle \gamma'\rangle$ 
is quasiregular.

\hfill

{\bf Proof:} Since $\gamma'$ is 
close to $\gamma\notin \tilde G_K$,
we may assume  that $\gamma'$ does not lie in the
compact subgroup $\tilde G_K \subset \tilde G_\C$.
Denote by $v'$ the vector $\log \gamma'$, and let
 $v''=\lambda_0 v'$
be the integer point on the line $\C v'$ which exists
by our assumptions. Since $v''\in \Lie(\tilde G_K)$, the 
coefficient $\lambda_0$ is not real.

Consider the holomorphic Lee flow 
$\psi'_\lambda= e^{a v'}$ on $\tilde M/\langle \gamma'\rangle$.
Clearly, $\psi'_a$ acts on $M$ trivially
for $a =1$ (because $e^{v'}=\gamma'$)
and $a=a_0$ (because $v''=a_0 v'$
is integer in $\Lie(\tilde G_K$). Therefore, the Lee flow is factorized 
through
the action of a compact Lie group
 $\C/\langle 1, a_0\rangle$,
and has compact orbits. We proved
that $\tilde M/\langle \gamma'\rangle$
is quasiregular. \endproof

%%%%%%%%%%%%%%%%%%%%%%%%%%%%%%%%%%%%%%%%%%%%%%%%

\section{Immersion theorem for general Vaisman manifolds}

%%%%%%%%%%%%%%%%%%%%%%%%%%%%%%%%%%%%%%%%%%%%%%%%

The main result of this paper is the following

\hfill

%%%%%%%%%%%%%%%%%%%%%%%%%%%%%%%%%%%%%%%%%%%%%%%%
\theorem \label{_imme_the_gene_Theorem_}
Let $M$ be a compact Vaisman manifold. Then
$M$ admits an immersion $M\stackrel \lambda \arrow H_\Lambda$
to a Hopf manifold. Moreover, $\lambda$ is 
compatible with the action of the Lee flow on $M$ and $H_\Lambda$,
and induces a finite covering $\lambda:\; M\arrow \lambda(M)$
{}from $M$ to its image.

\hfill

%%%%%%%%%%%%%%%%%%%%%%%%%%%%%%%%%%%%%%%%%%%%%%%%%%%%%%%%%%%%
\remark 
The converse statement is also true: given a compact 
submanifold $X \subset H_\Lambda$, $\dim X>0$, 
the manifold $X$ is LCK (which is clear from
the definitions) and Vaisman 
(see \cite{_Verbitsky:LCHK_}, Proposition 6.5).

\hfill

%%%%%%%%%%%%%%%%%%%%%%%%%%%%%%%%%%%%
\remark
\ref{_imme_for_qr_Vaisman_Theorem_} (Immersion Theorem
for quasiregular  Vaisman manifolds)
is a special case of \ref{_imme_the_gene_Theorem_}.

\hfill

\noindent{\bf Proof of \ref{_imme_the_gene_Theorem_}:} 
Let $\tilde M$ be the weight covering of $M$,
$\Gamma\cong Z$ its deck transform group,
$M\cong \tilde M/\Gamma$, and let $\gamma$ be
a generator of $\Gamma$.  By \ref{_qr_defo_exi_Proposition_}, 
there exists $\gamma'\in \tilde G_\C \subset \Aut(\tilde M, M)$ 
sufficiently close to $\gamma$, such that the quotient
$M':= \tilde M/\langle \gamma'\rangle$ is a quasiregular 
Vaisman manifold. Applying the same argument we used
to construct \eqref{_tilde_M_immersion_CD_Equation_},
we find an immersion
\[ \lambda':\; \tilde M \arrow \C^n\backslash 0
\]
with the following properties.
\begin{description}
\item[(i)] $\lambda'$ is a finite covering of its image
\item[(ii)] $\lambda'$ is a map associated with the space 
$H^0(M', L_C^k)$ of all holomorphic functions on $\tilde M$
which satisfy the following condition
\begin{equation}\label{_H^0(L^k)_defini_via_mono_Equation_}
H^0(M', L_C^k)= \{\chi:\; \tilde M \arrow \C \ \ |\ \ 
  \chi_i \circ \gamma' = q^k \chi\}
\end{equation}
where $q>1$ is monodromy action of $\gamma'$ in 
the weight bundle on $M'$.
\end{description}
By construction, $\gamma$ and $\gamma'$ commute
(they both sit in the same commutative group $\tilde G_\C$;
see Subsection \ref{_qr_Vaisman_defo_Subsection_}).
Therefore, the action of $\gamma$ on $\tilde M$
preserves the space \eqref{_H^0(L^k)_defini_via_mono_Equation_}.
This gives a commutative square

\begin{equation}\label{_imme_tilde_M_CD_Equation_}
\begin{CD}
\tilde M @>{\gamma}>> \tilde M \\
@V{\lambda'}VV  @VV{\lambda'}V \\
H^0(M', L_\C^k) @>{\lambda'(\gamma)}>> H^0(M', L_\C^k).
\end{CD}
\end{equation}
If $\gamma'$ is chosen sufficiently close to
$\gamma$, the action of $\lambda'(\gamma)$
on $H^0(M', L_\C^k)$ will be sufficiently close
to $q^k$, which can be chosen arbitrarily big.
Therefore, we may assume that 
all eigenvalues of $\lambda'(\gamma)$
are $>1$. This implies that
\[ H:= \left(H^0(M', L_\C^k)\backslash 0\right)/\lambda'(\gamma)
\]
is a well-defined complex manifold. From 
\eqref{_imme_tilde_M_CD_Equation_} we obtain
an immersion
\[\tilde M /\langle \gamma\rangle \arrow 
   \left(H^0(M', L_\C^k)\backslash 0\right)/\lambda'(\gamma),
\]
that is, an immersion of $M$ to $H$ which is 
a finite covering of its image.

To finish the proof of \ref{_imme_the_gene_Theorem_},
it remains to show that $H$ is a Hopf manifold, that is,
to show that the operator $\lambda'(\gamma)$ is semisimple
(\ref{_Hopf_mfld_Definition_}). 
This is implied by the following claim.

\hfill

%%%%%%%%%%%%%%%%%%%%%%%%%%%%%%%%%%%%%%%%%%%%%%%%%%%%%%%%%%%%%%%%%%%%%%%%
\claim\label{_gamma_semisimple_on_eigenva_of_gamma'_Claim_}
Let $M$ be a compact Vaisman manifold, $\tilde M$ its 
weight covering, and $\tilde G_\C$ the complex 
Lie group acting on $\tilde M$ as in
Subsection \ref{_tilde_G_C_Subsection_}. 
Take arbitrary elements $\gamma, \gamma' \in \tilde G_\C$,
and denote by $V_{\gamma'}$ the vector space of all
holomorphic functions $f:\; \tilde M \arrow \C$
satisfying 
\[ f\circ \gamma' = q^k f.
\]
  Assume that $V_{\gamma'}$ is finite-dimensional.
Then the natural action of $\gamma$ on $V_{\gamma'}$ is
semisimple. 

\hfill

\noindent{\bf Proof:} Consider the natural action of 
$\tilde G_\C$  on $\tilde M$. This action
is holomorphic and commutes with $\gamma'$, hence
$\tilde G_\C$ acts on $V_{\gamma'}$. Now, $\gamma\in \tilde G_\C$
is an element of a group $\tilde G_\C \cong (\C^*)^n$ acting on a 
finite-dimensional vector space.
All elements of
a subgroup $G\subset GL(n, \C)$, $G\cong (\C^*)^n$  
are semisimple, as an elementary group-theoretic
argument shows. Therefore, $\gamma$ is semisimple.
\ref{_gamma_semisimple_on_eigenva_of_gamma'_Claim_}
is proven and $H$ is in fact a Hopf manifold $H_\Lambda$. 
This proves \ref{_imme_the_gene_Theorem_}. \endproof

\hfill

\remark\label{_induced_structure_} It follows from \cite[Th. 
5.1]{_Vaisman:Dedicata_}, \cite[Th. 3.2]{_Tsukada_} and
\cite[Pr. 6.5]{_Verbitsky:LCHK_} that the metric induced on
$M$ by the  above holomorphic immersion is a Vaisman metric. But it
doesn't necessarily coincide with the initial one. Hence the immersion
is not isometric.

%%%%%%%%%%%%%%%%%%%%%%%%%%%%%%%%%%%%%%%%%%%%%%%%%%%%%%%%%%%%

\section{Applications to Sasakian geometry}

%%%%%%%%%%%%%%%%%%%%%%%%%%%%%%%%%%%%%%%%%%%%%%%%%%%%%%%%%%%%

As a by-product of the immersion theorem for compact Vaisman
manifolds, we derive a similar result for compact Sasakian
manifolds, the model space now being the odd sphere equipped with a
generally non-round metric and with a Sasakian structure obtained by
deforming the standard one by means of an $S^1$ action (see
\cite{_Kamishima_Ornea_}). We briefly recall this construction:

Let $S^{2n-2}$ be the unit sphere of $\C^n$ endowed with its standard
round metric and contact structure $\eta_0=\sum
(x_jdy_j-y_jdx_j)$. Let also $J$ denote the almost complex structure
of the contact distribution.  Deform $\eta_0$ by means of an action of $S^1$ 
by
setting $\eta_A=\frac{1}{\sum a_j|z_j|^2}\eta_0$, for $0<a_1\leq
a_2\cdots\leq a_n$.  Its Reeb field is $R_A=\sum a_j(x_j\partial
y_j-y_j\partial x_j)$. Clearly, $\eta_0$ and $\eta_A$ underly the same
contact structure. Finally, define the metric $g_A$ as follows:
\begin{itemize}
\item $g_A(X,Y)=d\eta_A(IX,Y)$ on the
  contact distribution~;
\item $R_A$ is normal to the contact distribution and has unit length.
\end{itemize}
It can be seen  that $S^{2n-1}_A:=(S^{2n-1}, g_A)$ is a Sasakian
manifold. If all the $a_j$ are equal, this structure corresponds to a
homothetic transformation along the contact distribution (a
D-homothetic transformation in the terms of S. Tanno \cite{_Tanno_}).

The Hopf manifold $H_\Lambda$, with $\lambda_j=e ^{-a_j}$, is obtained by 
identifying
$S^{2n-1}_A\times \R$ with $\C^n\setminus 0$ by means of the
diffeomorphism $(z_j, t)\mapsto (e ^{-a_jt}z_j)$. In particular, the
Lee form of  $H_\Lambda$ can now be identified with $-dt$, hence the
Lee flow corresponds to the action of the $S^1$ factor.

\hfill

Let now  $(X,h)$ be a compact Sasakian manifold. Then the trivial
bundle $M=X\times S^1$ has a Vaisman structure $(I,g)$ which is
quasiregular if the Reeb field $\xi$ of $X$ (which, on $M$, corresponds
to $I\theta ^\sharp$), determines a regular foliation.  From
\ref{_imme_the_gene_Theorem_}, $M$ is immersed in a Hopf manifold
$H_\Lambda$ and the immersion commutes with the action of the
respective Lee flows. As the Lee field of $H_\Lambda$ is induced from
the action of the $S^1$ factor, the immersion descends to an immersion
of Sasakian manifolds $f:X\arrow S_A ^{2n-1}$. If the Reeb field if $X$
is regular, then $M$ is quasiregular and $H_\Lambda$ is a diagonal
Hopf manifold, hence all $a_j$ are equal.

Note that the Reeb field of $X$ corresponds on the Vaisman manifold
$M$ to the anti-Lee field $I\theta ^\sharp$. As the holomorphic
immersion of Vaisman manifolds preserves the Lee field, the induced
immersion at the Sasakian level preserves the Reeb fields (still not
being an isometry).

On the other hand, recall that, by construction, the restriction of the
complex structure $I$ to the contact distribution $D$ of $X$ coincides
with the natural CR-structure of $X$ (see
\cite{_Dragomir_Ornea_}). Hence, any $X\in \Gamma(D)$ is of
the form $X=IY$ with $Y\in \Gamma(D)$ and, as $df$ commutes with $I$, we
have $g_A(dfX, R_A)=g_A(dfIY, R_A)=g_A(I_0dfY, R_A)=0,$ because
$I_0R_A=0$ by on $S_A^{2n-1}$. Thus $dfX$ belongs to the contact
distribution of $S_A^{2n-1}$.

%coincides with
%the $(1,1)$-endomorphism $\phi:=\nabla^h \xi$

Summing up we proved:

\hfill

%%%%%%%%%%%%%%%%%%%%%%%%%%%%%%%%%%%%%%%%%%%%%%%%%%%%%%%%%%%%
\theorem\label{_sas_imm_} Let $X$ be a compact Sasakian manifold.
Then $X$ admits an immersion into a Sasakian odd sphere  $S_A
^{2n-1}$ as above. This immersion is compatible with the
contact structure on $X$, $S_A^{2n-1}$, and 
preserves the Reeb field.

\hfill

The existence of an effective $U(1)$ action on a compact Riemannian
manifold imposes severe geometric and topological restrictions. From 
\ref{_qr_defo_exi_Proposition_} we derive that any
compact Sasakian manifold admits such an action. More precisely:

\hfill

%%%%%%%%%%%%%%%%%%%%%%%%%%%%%%%%%%%%%%%%%%%%%%%%%%%%%%%%%%%%
\proposition
Any compact Sasakian manifold $X$ has a deformation which is
isomorphic to a circle bundle over a projective orbifold.
In particular, any Sa\-sa\-kian manifold
admits an effective circle action.

\hfill

We note that this result is in fact true in a more general
context, for $K$-contact manifolds, as proved in
\cite{_Rukimbira_}.

Recall that a compact manifold which has an effective
circle action necessarily has zero minimal volume (cf. \cite{_Gromov_};
we are thankful to Gilles Courtois for explaining us this point). On the
other hand, a compact manifold which admits a metric of strictly
negative sectional curvature must have non-zero minimal volume, {\em
loc. cit.} We obtain the following corollary, previously proved
in \cite[Pr. 2.7]{_Boyer_Galicki:proc_} for $K$-contact manifolds and,
using harmonic maps, for Sasakian manifolds, in \cite{_Petit_}:

\hfill

\corollary A compact Sasakian manifold cannot bear any Riemannian metric of 
strictly negative
sectional curvature.

\hfill

{\small

\hfill

\noindent {\sc Liviu Ornea\\
University of Bucharest, Faculty of Mathematics, 14
Academiei str., 70109 Bucharest, Romania.}\\
\tt lornea@imar.ro

\hfill

\noindent {\sc Misha Verbitsky\\
University of Glasgow, Department of Mathematics, 15
  University Gardens, Glasgow, Scotland.}\\
\tt verbit@maths.gla.ac.uk, \ \  verbit@mccme.ru 
}% end of small

\end{document}